\documentclass[12pt,letterpaper]{article}
\usepackage{amsthm,amssymb}
\usepackage{amsfonts}
\makeatletter
\usepackage[dvips]{color}
\usepackage{colordvi,multicol}

\newtheorem{thm}{\textbf Theorem}[section]
\newtheorem{lem}[thm]{\textbf Lemma}
\newtheorem{rem}[thm]{\textbf Remark}

\newcommand{\be}{\begin{eqnarray*}}
\newcommand{\ee}{\end{eqnarray*}}

\begin{document}

\title{\bf Probabilistic representations of solutions of elliptic boundary value problem and non-symmetric semigroups}
  \author{
Chuan-Zhong Chen\\
School of Mathematics and Statistics\\
Hainan Normal University\\
Haikou 571158, China\\
czchen@hainnu.edu.cn\\
\\
Wei Sun
            \\ Department of Mathematics and Statistics\\
             Concordia University\\
             Montreal H3G 1M8, Canada\\
             wei.sun@concordia.ca\\
\\
Jing Zhang
 \\ Department of Mathematics and Statistics\\
             Concordia University\\
             Montreal H3G 1M8, Canada\\
             waangel520@gmail.com}
    \date{}
\maketitle

\begin{abstract}
\noindent In this paper, we use a probabilistic approach to show
that there exists a unique, bounded continuous solution to the
Dirichlet boundary value problem for a general class of second
order non-symmetric elliptic operators $L$ with singular
coefficients, which does not necessarily have the maximum
principle. The theory of Dirichlet forms and heat kernel estimates
play a crucial role in our approach. A probabilistic
representation of the non-symmetric semigroup $\{T_t\}_{t\ge 0}$
generated by $L$ is also given.

\end{abstract}

\section{Introduction and the Main Theorem}

In this paper, we will use probabilistic methods to study the
Dirichlet boundary value problem for second order elliptic
differential operators:
\begin{eqnarray}\label{111}
\left\{
  \begin{array}{ll}
   L u=0 \ \ \ \ {\rm in}\ \ D \\
    u=f \ \ \ \ \ \ {\rm on}\ \  \partial D,
  \end{array}
\right.
\end{eqnarray}
where $D$ is a bounded connected open subset of $\mathbb{R}^d$.
The operator $L$ is given by
\begin{equation}\label{12}
Lu=\frac{1}{2}\sum_{i,j=1}^d\frac{\partial}{\partial
x_j}\left(a_{ij}(x)\frac{\partial u}{\partial
x_i}\right)+\sum_{i=1}^db_i(x)\frac{\partial u}{\partial x_i}+
(c(x)-{\rm div}\,\hat b(x))u,
\end{equation}
where $A(x)=(a_{ij}(x))_{i,j=1}^d$ is a Borel measurable, (not
necessarily symmetric) matrix-valued function  on $D$ satisfying
\begin{equation}\label{13}
\lambda|\xi|^2\le \sum_{i,j=1}^da_{ij}(x)\xi_i\xi_j\ \ {\rm for\
any}\ \xi=(\xi_i)_{i=1}^d\in \mathbb{R}^d, x\in D
\end{equation}
and
\begin{equation}\label{14}
|a_{ij}(x)|\le \frac{1}{\lambda}\ \ {\rm for\ any}\ x \in D,\ 1\le
i,j\le d
\end{equation}
for some constant $0<\lambda\leq1$; $b=(b_1,\dots,b_d)^*$ and
$\hat{b}=(\hat{b}_1,\dots,\hat{b}_d)^*$ are Borel measurable
$\mathbb{R}^d$-valued functions on $D$ and $c$ is a Borel
measurable function on $D$ satisfying $|b|^2\in L^{p\vee
1}(D;dx)$, $|\hat b|^2\in L^{p\vee 1}(D;dx)$ and $c\in L^{p\vee
1}(D;dx)$ for some constant $p>d/2$. Hereafter we use $^*$ to
denote the transpose of a
 vector or matrix, and use $|\cdot|$ and $\langle\cdot,\cdot\rangle$ to denote
respectively the standard norm and inner product of the Euclidean
space $\mathbb{R}^d$.

In (\ref{111}), $Lu=0$ in $D$ is understood in the distributional
sense:
$$
u\in H^{1,2}(D)\ {\rm and}\ {\cal E}(u,\phi)=0\ {\rm for\ every}\
\phi\in C^{\infty}_0(D),
$$
where $H^{1,2}(D)$ is the Sobolev space on $D$ with norm
$$\|f\|_{H^{1,2}}:=\left(\int_{D}|\nabla
f(x)|^2dx+\int_{D}|f(x)|^2dx\right)^{1/2},$$ $C^{\infty}_0(D)$ is
the space of infinitely differentiable functions with compact
support in $D$, and $({\cal E},D(\cal{E}))$ is the bilinear form
associated with $L$:
\begin{eqnarray}\label{EEE}
{\cal
E}(u,v)&=&\frac{1}{2}\sum_{i,j=1}^d\int_{D}a_{ij}(x)\frac{\partial
u}{\partial x_i}\frac{\partial v}{\partial x_j}dx
-\sum_{i=1}^d\int_{D}b_i(x)\frac{\partial u}{\partial x_i}v(x)dx\nonumber\\
& &-\sum_{i=1}^d\int_{D}\hat{b}_i(x)\frac{\partial (uv)}{\partial x_i}dx-\int_{D}c(x)u(x)v(x)dx,\\
D({\cal E})&=&H^{1,2}_0(D)\nonumber
\end{eqnarray}
with $H^{1,2}_0(D)$ being the completion of $C_0^{\infty}(D)$ with
respect to the Sobolev-norm $\|\cdot\|_{H^{1,2}}$.
 By setting $a=I,\ b=0,\
\hat b=0$ and $c=0$ off $D$, we may assume that the operator $L$
is defined on $\mathbb{R}^d$.

Using probabilistic approaches to solve boundary value problems
has a long history. The pioneering work goes back to Kakutani
\cite{Kakutani}, who used Brownian motion to represent the
solution of the classical Dirichlet boundary value problem with
operator  $L=\Delta$, the Laplacian operator. If $\hat b=0$ and
$c\le 0$, then the solution $u$ to problem (\ref{111}) is given by
the famous Feynman-Kac formula
$$
u(x)=E_x\left[e^{\int_0^{\tau_D}c(X_s)ds}f(X_{\tau_D})\right], \ \
x\in D,
$$
where $X=(X_t)_{t\ge 0}$ is the diffusion process associated with
the generator $L^b$ given by
\begin{equation}\label{Lbb}
L^bu=\frac{1}{2}\sum_{i,j=1}^d\frac{\partial}{\partial
x_j}\left(a_{ij}(x)\frac{\partial u}{\partial
x_i}\right)+\sum_{i=1}^db_i(x)\frac{\partial u}{\partial x_i},
\end{equation}
and $\tau_D$ is the first exit time of $X$ from $D$. We refer the
readers to \cite{Chen} for the general results obtained in this
case.

When $\hat b\not=0$ and $A$ is symmetric, Chen and Zhang
\cite{Zhang3} used the time reversal of symmetric Markov processes
to give an explicit probabilistic representation of the solution
to problem (\ref{111}). (Note that the operator $L$ given by
(\ref{12}) is the same as that used in \cite{Zhang3} if we replace
$b$ with $b-\hat b$ in (\ref{12}).) We should point out that the
${\rm div}\,\hat b$ in (\ref{12}) is just a formal writing since
the vector field $\hat b$ is merely measurable hence its
divergence exists only in the distributional sense. In the
remarkable paper \cite{Zhang3}, Chen and Zhang proved that there
exists a unique, bounded continuous weak solution to  problem
(\ref{111}) without the Markov assumption
\begin{equation}\label{Markov}
c-{\rm div}\,\hat b\le 0\ \ {\rm in}\ \mathbb{R}^d,
\end{equation}
i.e.,
$\int_{\mathbb{R}^d}c(x)\phi(x)dx+\sum_{i=1}^d\int_{\mathbb{R}^d}\hat{b}_i(x)\frac{\partial
\phi}{\partial x_i}dx\le 0$ for any nonnegative $\phi\in
C^{\infty}_0(\mathbb{R}^d)$. The novelty of \cite{Zhang3} is to
tackle the lower-order term ${\rm div}\,\hat b$ through combining
the time-reversal of a Girsanov transform from the random time
$\tau_D$ with a certain $h$-transform. In \cite{Zhang3}, Chen and Zhang used
essentially the following result due to Meyers \cite{Meyers}:

{\it For every $x_0\in \mathbb{R}^d$, $R>0$ and $p>d$, there is a
constant $\varepsilon\in (0,1)$, depending only on $d$, $R$ and
$p$, such that if
\begin{equation}\label{3.1}
(1-\varepsilon)I_{d\times d}\le A(x)\le I_{d\times d}\ \ {\rm for\
a.e.}\ x\in B_R:=B(x_0,R),
\end{equation}
then \begin{equation}\label{VMO}\frac{1}{2}\nabla(A\nabla u)={\rm
div}\,f\end{equation} in $B_R$ has a
unique weak solution in $H^{1,p}_0(B_R)$ for every
$f=(f_1,\dots,f_d)\in L^p(B_R;dx)$. Moreover, there is a constant
$c>0$ independent of $f$ such that $$ \|\nabla
u\|_{L^p(B_R;dx)}\le c\|f\|_{L^p(B_R;dx)}.
$$}

To apply Meyers's result,  the diffusion matrix $A$ is assumed to
satisfy Condition (\ref{3.1}) in \cite{Zhang3}  (see
\cite[Theorems 3.3 and 4.5]{Zhang3}). If Condition (\ref{3.1}) is
replaced with other conditions which guarantee that $\nabla u$ in
(\ref{VMO}) belongs to some $L^p$ space for $p>d$, e.g. the
condition that $A$ is in the class VMO and $\partial D\in C^{1,1}$
(see \cite{Faz}), then Chen and Zhang's approach still apply. We
thank Professors Z.Q. Chen and T.S. Zhang for pointing out this to
us.

In general, it is possible that $f\in L^p$ while $\nabla u\notin
L^p$ (see \cite{Meyers} for an example). For this case, we cannot
use the $h$-transform method to tackle the lower-order term ${\rm
div}\,\hat b$ even when $A$ is symmetric. In this paper, we will
show that there exists a unique, bounded continuous solution to
problem (\ref{111}) without additional condition on $A$ such as
Condition (\ref{3.1}), the VMO condition or the symmetry of $A$,
and without the Markovian assumption (\ref{Markov}). Instead of
using Meyers's $L^p$-estimate as in \cite{Zhang3}, we will make
use of Aronson's heat kernel estimates \cite{Aron1,Aron}.

In the sequel, we  let $X=((X_t)_{t\ge 0},(P_x)_{x\in
\mathbb{R}^d})$ be the Markov process associated with the
following (non-symmetric) Dirichlet form
\begin{eqnarray}\label{EEF}
{\cal E}^0(u,v)&=&\frac{1}{2}\int_{\mathbb{R}^d}\sum_{i,j=1}^da_{ij}(x)\frac{\partial u}{\partial x_i}\frac{\partial v}{\partial x_j}dx,\\
D({\cal E}^0)&=&H^{1,2}(\mathbb{R}^d).\nonumber
\end{eqnarray}
It is well-known that $X$ is a conservative Feller process on
$\mathbb{R}^d$ that has continuous transition density function
which admits a two-sided Aronson's heat kernel estimate. Let
$\{{\cal F}_t,t\ge 0\}$ be the minimal augmented filtration
generated by $X$. By Fukushima's decomposition (cf. \cite[Theorem 5.1.8]{oshima}), we have
$$
X_t=x+M_t+N_t,
$$
where $M_t=(M_t^1,\dots,M_t^d)^*$ is a martingale additive
functional of $X$ with quadratic co-variation
$$
\langle M^i,M^j\rangle_t=\int_0^t \tilde a_{ij}(X_s)ds
$$
and $N_t=(N_t^1,\dots,N_t^d)^*$ is a continuous additive
functional of $X$ locally of zero quadratic variation. Hereafter
$\tilde A=(\tilde a_{ij})_{i,j=1}^d$ denotes the symmetrization of
$A$, i.e., $\tilde A:=1/2(A+A^*)$.

For any vector field $\xi\in L^{2}(\mathbb{R}^d;dx)$,
 there exists a unique function $\xi^H\in H^{1,2}(\mathbb{R}^d)$
 such that
 $$\int_{\mathbb{R}^d}\langle \xi,\nabla h\rangle dx=-{\cal E}_1^0(\xi^H,h),\ \ \forall h\in C^{\infty}_0(\mathbb{R}^d)$$
 (see Lemma \ref{lem111} below). Hereafter ${\cal E}^0_{\gamma}(u,v):={\cal E}^0(u,v)+\gamma\int uvdx$ for any
 $u,v\in D({\cal E}^0)$ and any
 constant $\gamma$.
 We have Fukushima's decomposition:
$$
\widetilde{{\xi}^H}(X_t)-\widetilde{{\xi}^H}(X_0)=M^{{\xi}^H}_t+N^{{\xi}^H}_t,
$$
where $\widetilde{{\xi}^H}$ is a quasi-continuous version of
${\xi}^H$. To simplify notation, in the sequel we take $w$ to be
its quasi-continuous version $\tilde w$ whenever such a version
exists. As in \cite{Fu94,MR92}, we use the term ``quasi-everywhere" (abbreviated ``q.e.") to mean ``except on an exceptional set".

Now we can state the main theorem of this paper.

\begin{thm}\label{main1} Let $d\geq1$, $D$ be a bounded Lipschitz domain in $\mathbb{R}^d$ and $p>d/2$.
Suppose that

(i) $A$ satisfies (\ref{13}) and (\ref{14}).

(ii) $|b|^2\in L^{p\vee 1}(D;dx)$ and $|\hat b|^2\in L^{p\vee
1}(D;dx)$.

(iii) $c\in L^{p\vee 1}(D;dx)$ and $c-{\rm div}\,\hat b\leq g$ for
some nonnegative function $g\in L^{p\vee 1}(D;dx)$ in the
distributional sense.

Then, there exists a constant $M>0$ such that whenever
$\|g\|_{L^{p\vee 1}}\leq M$, for any $f\in C(\partial D)$, there
exists a unique weak solution $u$ to $Lu=0$ in $D$ that is
continuous on $\overline{D}$ with $u=f$ on $\partial D$. Moreover,
the solution $u$ admits the following representation: for q.e.
$x\in D$,
\begin{eqnarray}\label{15}
u(x)&=&E_x\left[\exp\left(\int_0^{\tau_D}(\tilde
a^{-1}b)^*(X_s)dM_s-\frac{1}{2}\int_0^{\tau_D}
b^*\tilde a^{-1}b(X_s)ds\right.\right.\nonumber\\
& &\ \ \ \ \left.\left. +\int_0^{\tau_D} c(X_s)ds+N^{{\hat b}^H}_{{\tau_D}}-\int_0^{\tau_D} {\hat b}^H(X_s)ds\right)f(X_{{\tau_D}})\right].
\end{eqnarray}
\end{thm}

We will give the proof of Theorem \ref{main1} in Section 2, which
consists of three subsections. In Subsection 2.1, we prove the
existence of the weak solution and give its probabilistic
representation (\ref{15}). In Subsection 2.2, we prove the
continuity of the weak solution.  In Subsection 2.3,  we prove the
uniqueness of the continuous weak solutions. The recently
developed Nakao integral for non-symmetric Dirichlet forms (cf.
\cite{Alexander 2013} and \cite{CMS}) will be used in the proof of
the uniqueness.

In Section 3, we use some techniques of Section 2 to give a
probabilistic representation of the non-symmetric semigroup
$\{T_t\}_{t\ge 0}$ generated by $L$ that is defined by (\ref{12}).
The obtained result (see Theorem \ref{main} below) generalizes the
corresponding result of \cite{Zhang1} from the case of symmetric
diffusion matrix $A$ to the non-symmetric case.

\section{Proof of Theorem \ref{main1}}\setcounter{equation}{0}

\subsection{Proof of the existence of weak solution}\setcounter{equation}{0}

We first generalize \cite[Theorem 1.1]{Chen} from the case of
symmetric diffusion matrix $A$ to the non-symmetric case. Define
$$
L^1u=\frac{1}{2}\sum_{i,j=1}^d\frac{\partial}{\partial
x_j}\left(a_{ij}(x)\frac{\partial u}{\partial
x_i}\right)+\sum_{i=1}^db_i(x)\frac{\partial u}{\partial x_i}+
c(x)u.
$$
\begin{lem}\label{05}
Suppose that $D$ is a bounded domain in $\mathbb{R}^d$, $c\le 0$
and $f\in C(\partial D)$. Then
\begin{eqnarray*}
u(x)&=&E_x\left[\exp\left(\int_0^{\tau_D}(\tilde
a^{-1}b)^*(X_s)dM_s-\frac{1}{2}\int_0^{\tau_D}
b^*\tilde a^{-1}b(X_s)ds\right.\right.\nonumber\\
& &\ \ \ \ \left.\left. +\int_0^{\tau_D}
c(X_s)ds\right)f(X_{{\tau_D}})\right]
\end{eqnarray*}
is the unique weak solution of $L^1u=0$ which is continuous in $D$
and
$$\lim_{x\rightarrow y,x\in D}u(x)=f(y)$$
for $y\in \partial D$ which is regular for the Laplace operator $(\frac{1}{2}\triangle, D)$.
\end{lem}

\noindent {\bf Proof.} The proof of Lemma \ref{05} is similar to
that of \cite[Theorem 1.1]{Chen}. We only point out below the main
differences in the argument between the symmetric and the
non-symmetric cases.

Denote by $X^0$ the part of the process $X$ on $D$, that is, $X^0$
is obtained by killing the sample paths of $X$ upon leaving $D$.
By \cite{Aron1, Aron}, the transition
density function $p_0(t,x,y)$ of $X^0$ has the upbound estimate
\begin{equation}\label{30}
    p_0(t,x,y)\leq
    \frac{\vartheta}{t^{d/2}}e^{-\frac{|x-y|^2}{\vartheta t}},\ \ (t,x,y)\in (0,\infty)\times D\times D,
\end{equation}
for some constant $\vartheta>0$.

We define
$$
L^0u=\frac{1}{2}\sum_{i,j=1}^d\frac{\partial}{\partial
x_j}\left(a_{ij}(x)\frac{\partial u}{\partial x_i}\right).
$$
Let $D_1$ be a bounded subdomain of $D$ and $f_1\in
H_{0}^{1,2}(D)$. By \cite{Trudinger}, there exists a unique weak
solution of $L^0u=0$ in $D_1$ such that $u-f_1|_{D_1}\in
H_{0}^{1,2}(D_1)$. Further, by the famous theorem of Littman,
Stampacchia and Weinberger, which holds also for the non-symmetric
case (cf. e.g. \cite{Kenig}), we can prove the analog of
\cite[Theorem 2.1]{Chen} with the non-symmetric $A$. By virtue of
the Harnack inequality for parabolic equations (cf. \cite{Moser}
and \cite{Lierl}), we can prove that \cite[Lemma 2.2]{Chen} and
hence \cite[Corollary 2.3 and Theorem 2.4]{Chen} hold for the
non-symmetric case.

Finally, we would like to point out that the exponential
martingale  $M_t$ introduced in \cite[(3.4)]{Chen} needs to be
replaced with
\begin{equation}\label{U}
U_t:=\exp\left(\int_0^t(\tilde
a^{-1}b)^*(X_s)dM_s-\frac{1}{2}\int_0^tb^*\tilde
a^{-1}b(X_s)ds\right),\ \ t\geq0
\end{equation}
for our non-symmetric case. \hfill\fbox

\begin{lem}\label{lem111}
(i) For any vector field $\xi\in L^2(\mathbb{R}^d;dx)$, there
exists a unique function $\xi^H\in H^{1,2}(\mathbb{R}^d)$ such
that
\begin{equation}\label{2002}
\int_{\mathbb{R}^d}\langle \xi,\nabla h\rangle dx=-{\cal
E}_1^0(\xi^H,h),\ \ \forall  h\in H^{1,2}(\mathbb{R}^d).
\end{equation}

(ii) If $\xi_n$ converges to $\xi$ in $L^2(\mathbb{R}^d;dx)$ as
$n\rightarrow\infty$, then $\xi^H_n$ converges to $\xi^H$ in
$H^{1,2}(\mathbb{R}^d)$ as $n\rightarrow\infty$.

 (iii) For $\xi\in
C_0^{\infty}(\mathbb{R}^d)$,
\begin{equation}\label{2004}
  -\int_0^t {\rm div}\,\xi(X_s)ds=N^{{\xi}^H}_t-\int_0^t {\xi}^H(X_s)ds,\ \ t\ge 0.
\end{equation}

\end{lem}

\noindent {\bf Proof.}  (i) Let $\xi\in L^2(\mathbb{R}^d;dx)$. We
define the map $\eta: h\in H^{1,2}(\mathbb{R}^d)\mapsto
\int_{\mathbb{R}^d}\langle \xi,\nabla h\rangle dx$. By the Riesz
representation theorem, there exists a unique $\xi^0\in
H^{1,2}(\mathbb{R}^d)$ such that
\begin{equation}\label{may14}
\eta(h)=\tilde{\cal E}_1^0(\xi^0,h),\ \ \forall h\in
H^{1,2}(\mathbb{R}^d),
\end{equation}
where $({\tilde{\cal E}}^0, D({\cal E}^0))$ denotes the symmetric
part of the Dirichlet form $({\cal E}^0, D({\cal E}^0))$. Thus, by
\cite[Lemma 2.1]{CMS}, there exists a unique $\xi^H\in D({\cal
E}^0)= H^{1,2}(\mathbb{R}^d)$ such that
\begin{equation}\label{lemma222}
\tilde{\cal E}_1^0(\xi^0,h)=-{\cal E}_1^0(\xi^H,h),\ \ \forall
h\in H^{1,2}(\mathbb{R}^d).
\end{equation}

(ii) Suppose $\xi_n$ converges to $\xi$ in $L^2(\mathbb{R}^d;dx)$
as $n\rightarrow\infty$. By (\ref{may14}), we get
\begin{eqnarray}\label{chen}
\|\xi^0_n- \xi^0\|_{\tilde{{\cal E}}^0_1}&=& \sup_{\|h\|_{\tilde{\cal{E}}^0_1}=1} \tilde{\cal{E}}^0_1(\xi^0_n- \xi^0, h)\nonumber\\
& =& \sup_{\|h\|_{\tilde{\cal{E}}^0_1}=1}\int_{\mathbb{R}^d} \langle \xi_n- \xi, \nabla h\rangle dx \nonumber\\
& \leq & \| \xi_n- \xi
\|_{L^2}\sup_{\|h\|_{\tilde{\cal{E}}^0_1}=1}\|h
\|_{H^{1,2}}\nonumber\\
&\rightarrow&0\ {\rm as}\ n\rightarrow\infty.
\end{eqnarray}
Further, by (\ref{lemma222}), we get
\begin{eqnarray}\label{may1}
{\cal E}_1^0(\xi_n^H-\xi^H, \xi_n^H)&=&{\cal E}_1^0(\xi_n^H, \xi_n^H)- {\cal E}_1^0(\xi^H, \xi_n^H) \nonumber\\
& = & -\tilde{\cal E}_1^0(\xi_n^0, \xi_n^H) + \tilde{\cal E}_1^0(\xi^0, \xi_n^H)\nonumber\\
& = & \tilde{\cal{E}}_1^0(\xi^0-\xi_n^0, \xi_n^H)\nonumber\\
& \leq & \left[\tilde{\cal E}_1^0(\xi^0- \xi_n^0,
\xi^0-\xi_n^0)\right]^{1/2}\left[\tilde{\cal E}_1^0(\xi_n^H,
\xi_n^H)\right]^{1/2},
\end{eqnarray}
\begin{equation}\label{may2}\sup_{n\in\mathbb{N}} {\cal E}_1^0(\xi_n^H, \xi_n^H)\leq \sup_{n\in\mathbb{N}} \tilde{\cal E}_1^0(\xi_n^0, \xi_n^0)<\infty,\end{equation}
and \begin{eqnarray}\label{may3} \lim_{n\rightarrow\infty}{\cal
E}_1^0(\xi_n^H-\xi^H,\xi^H)&=&-\lim_{n\rightarrow\infty}\tilde{\cal
E}_1^0(\xi_n^0-\xi^0,\xi^H)\nonumber\\
&=&-\lim_{n\rightarrow\infty}\int_{\mathbb{R}^d}\langle \xi_n-\xi,\nabla \xi^H\rangle dx\nonumber\\
&=&0.
\end{eqnarray}
Therefore, we obtain by (\ref{chen})-(\ref{may3}) that
\begin{eqnarray*}
\lim_{n\rightarrow\infty}{\cal E}_1^0(\xi_n^H- \xi^H, \xi_n^H-
\xi^H) & =
& \lim_{n\rightarrow\infty}\left\{{\cal E}_1^0(\xi_n^H-\xi^H, \xi_n^H)- {\cal E}_1^0(\xi_n^H-\xi^H, \xi^H)\right\}\\
&=&0.
\end{eqnarray*}

(iii) Let $\xi\in C_0^{\infty}(\mathbb{R}^d)$. For any $h\in
H^{1,2}(\mathbb{R}^d)$, we have
\begin{eqnarray*}
\lim_{t\rightarrow 0}\frac{1}{t}E_{h\cdot dx}\left[-\int_0^t {\rm
div}\,\xi(X_s)ds\right]&=&
-\int_{\mathbb{R}^d}({\rm div}\,\xi)hdx\\
&=&\int_{\mathbb{R}^d}\langle \xi,\nabla h\rangle dx\\
&=&-{\cal E}_1^0(\xi^H,h)\\
&=&\lim_{t\rightarrow 0}\frac{1}{t}E_{h\cdot
dx}\left[N^{{\xi}^H}_t-\int_0^t {\xi}^H(X_s)ds\right].
\end{eqnarray*}
Therefore, (\ref{2004}) holds by \cite[Lemma 2.3]{CMS}.\hfill\fbox

\vskip 0.5cm \noindent {\bf Proof of the existence of weak
solution and its probabilistic representation.} \vskip 0.5cm

We define a family of measures $\{Q_x,x\in \mathbb{R}^d\}$ on
${\cal F}_\infty$ by
$$\left.\frac{d {Q}_x}{d{P}_x}\right|_{{\cal F}_t}=U_t,\ \ t\ge 0,$$
where $U_t$ is given by (\ref{U}). Then, under $\{Q_x,x\in
\mathbb{R}^d\}$, $X$ is a diffusion process on $\mathbb{R}^d$ with
the generator $L^b$ given by (\ref{Lbb}). Denote by $E^Q_x$ the
expectation with respect to the measure $Q_x$ for $x\in
\mathbb{R}^d$. \textbf{From now on till the end of this section,
we fix a constant $0<\theta<\frac{1}{2}$}. We will show below that
there exists a constant $M>0$ such that for any $w\in L^{p\vee 1}(\mathbb{R}^d;dx)$ with $\|w\|_{L^{p\vee 1}}\leq M$,
we have \begin{equation}\label{new333}\sup_{x\in
D}E_x^{Q}\left[\int_0^{\tau_D}|w|(X_s)ds\right]\leq \theta.
\end{equation}

 We only prove
(\ref{new333}) when $d\ge 3$. The cases that $d=1,2$ can be
considered similarly. Let $X^D$ be the part of the process $X$ on
$D$ under $\{Q_x\}$, that is, $X^D$ is obtained by killing the
sample paths of $X$ upon leaving $D$. Denote by $p(t,x,y)$ the
transition density function of $X^D$. By \cite[Theorem 9]{Aron},
for each $T>0$, there exist positive constants $\sigma^T_1$ and
$\sigma^T_2$ such that
$$
p(t,x,y)\leq \frac{\sigma^T_1}{t^{d/2}}e^{-\frac{\sigma^T_2|x-y|^2}{t}},\ \
(t,x,y)\in (0,T)\times D\times D.
$$
Similar to the proof of \cite[Lemma 6.1]{KS}, we  can show that
there exist positive constants $\sigma_1$ and $\sigma_2$
such that
\begin{eqnarray}\label{213}
p(t,x,y)\leq \frac{\sigma_1}{t^{d/2}}e^{-\frac{\sigma_2|x-y|^2}{t}},\ \
(t,x,y)\in (0,\infty)\times D\times D.
\end{eqnarray}
Denote by $G_D(x,y)$ the Green function of $X^D$. Then,
\begin{equation}\label{215}
    G_D(x,y)\leq \frac{\sigma_3}{|x-y|^{d-2}},\ \
(x,y)\in D \times D,
\end{equation}
for some positive constant $\sigma_3$.

Let $q>1$ satisfy $\frac{1}{p}+\frac{1}{q}=1$. Then $d-q(d-2)>0$.
We obtain by (\ref{215}) that
\begin{eqnarray*}
E_x^{Q}\left[\int_0^{\tau_D}|w|(X_s)ds\right]&=& \int_D G_D(x,y)|w|(y)dy\nonumber\\
&\leq&\int_D\frac{\sigma_3|w|(y)}{|x-y|^{d-2}}dy \nonumber\\
&\leq&\sigma_3\left(\int_D(|w|(y))^{p}dy\right)^{1/p}\left(\int_D|x-y|^{-q(d-2)}dy\right)^{1/q}\nonumber\\
&\leq&\sigma_3\|w\|_{L^{p}}\left(\int_0^\varsigma r^{d-q(d-2)-1}dr\right)^{1/q}\nonumber\\
&=&\frac{\sigma_3\varsigma^{d/q-(d-2)}}{[d-q(d-2)]^{1/q}}\|w\|_{L^{p}}.
\end{eqnarray*}
Hereafter $\varsigma$ denotes the diameter of $D$. Set
$$M:=\frac{\theta[d-q(d-2)]^{1/q}}{\sigma_3\varsigma^{d/q-(d-2)}}.$$
Then $\|w\|_{L^{p}}\leq M$ implies (\ref{new333}). Further, by
(\ref{new333}) and Khasminskii's inequality, we get
\begin{equation}\label{old334}
\sup_{x\in
D}E_x^{Q}\left[\exp\left(\int_0^{\tau_D}|w|(X_s)ds\right)\right]\le\frac{1}{1-\theta}.
\end{equation}

We define
$$
J(x)=\frac{1_{\{|x|<1\}}e^{-\frac{1}{1-|x|^2}}}{\int_{\{|y|<1\}}e^{-\frac{1}{1-|y|^2}}dy},
\ \ x\in \mathbb{R}^d.
$$
For $k\in \mathbb{N}$ and $x\in \mathbb{R}^d$, set
\begin{eqnarray*} J_k(x)&:=&k^dJ(kx),\\
\hat b_k(x)&:=&\int_{\mathbb{R}^d}\hat b(x-y)J_k(y)dy,\\
c_k(x)&:=&\int_{\mathbb{R}^d} c(x-y)J_k(y)dy,\\
g_k(x)&:=&\int_{\mathbb{R}^d} g(x-y)J_k(y)dy.
\end{eqnarray*}
We have
\begin{equation}\label{addsa11}
\hat b_k\rightarrow \hat b\ {\rm in}\ L^2(\mathbb{R}^d;dx)\ {\rm as}\ k\rightarrow\infty
\end{equation}
and
\begin{equation}\label{addsa}
c_k\rightarrow c\ {\rm in}\ L^1(\mathbb{R}^d;dx)\ {\rm as}\ k\rightarrow\infty.
\end{equation}
Suppose $\|g\|_{L^{p\vee 1}}\le M$. Since $c-{\rm div}\,\hat b\leq
g$ implies that $c_k-{\rm div}\,\hat b_k\leq g_k$ for $k\in
\mathbb{N}$, we obtain by (\ref{old334}) that
\begin{equation}\label{new334}
\sup_{k\in \mathbb{N}}\sup_{x\in
D}E_x^{Q}\left[\exp\left(\int_0^{\tau_D}(c_k-{\rm div}\,\hat
b_k)(X_s)ds\right)\right]\le\frac{1}{1-\theta}.
\end{equation}

Define for $t\ge 0$,
\begin{eqnarray}\label{311}
Z_t&:=&\exp\left(\int_0^t(\tilde
a^{-1}b)^*(X_s)dM_s-\frac{1}{2}\int_0^tb^*\tilde
a^{-1}b(X_s)ds\right.\nonumber\\
& &\ \ \left. +\int_0^{t} c(X_s)ds+N^{{\hat b}^H}_{{t}}-\int_0^{t}
{\hat b}^H(X_s)ds\right).
\end{eqnarray}
By (\ref{addsa11}) and Lemma \ref{lem111}(ii), we get
\begin{equation}\label{qe11}
\hat b_k^H\rightarrow \hat b^H\ {\rm in}\ H^{1,2}(\mathbb{R}^d)\
{\rm as}\ k\rightarrow\infty.
\end{equation}
Further, by \cite[Lemma 4.1.12 and Theorem
5.1.2]{oshima}, there exists a subsequence $\{k_l\}$ such that for q.e. $x\in \mathbb{R}^d$,
\begin{equation}\label{qe12}
P_x\left\{\lim_{l\rightarrow\infty}N^{\hat b_{k_l}^H}_t=N^{\hat
b^H}_t\ {\rm uniformly\ on\ any\ finite\ interval\ of}\
t\right\}=1.
\end{equation}
For simplicity, we still use $\{k\}$ to denote the subsequence $\{k_l\}$. By
(\ref{addsa})-(\ref{qe12}) and Fatou's lemma, we obtain  that
\begin{eqnarray}\label{ZtauD}
E_x[Z_{\tau_D}]&=& E_x^{Q}\left[\exp\left(\int_0^{\tau_D}c(X_s)ds+N_{\tau_D}^{\hat b^H}-\int_0^{\tau_D}\hat b^H(X_s)ds\right)\right]\nonumber \\
&\leq& \liminf_{k\rightarrow\infty}
E_x^{Q}\left[\exp\left(\int_0^{\tau_D}c_k(X_s)ds+N_{\tau_D}^{\hat
b_k^H}-
\int_0^{\tau_D}\hat b_k^H(X_s)ds\right)\right]\nonumber \\
&=& \liminf_{k\rightarrow\infty} E_x^{Q}\left[\exp\left(\int_0^{\tau_D}(c_k-{\rm div}\,\hat b_k)(X_s)ds\right)\right] \nonumber\\
&\le& \frac{1}{1-\theta},\ \ {\rm for\ q.e.}\ x\in D.
\end{eqnarray}

For $k\in \mathbb{N}$, we define
$$
L_ku=\frac{1}{2}\sum_{i,j=1}^d\frac{\partial}{\partial
x_j}\left(a_{ij}(x)\frac{\partial u}{\partial
x_i}\right)+\sum_{i=1}^db_i(x)\frac{\partial u}{\partial x_i}+
(c_k(x)-{\rm div}\,\hat b_k(x))u.
$$
The bilinear form  $({\cal E}_k,D({\cal E}_k))$ associated with
$L_k$ is
\begin{eqnarray*}
{\cal
E}_k(u,v)&=&\frac{1}{2}\sum_{i,j=1}^d\int_{D}a_{ij}(x)\frac{\partial
u}{\partial x_i}\frac{\partial v}{\partial x_j}dx
-\sum_{i=1}^d\int_{D}b_i(x)\frac{\partial u}{\partial x_i}v(x)dx\\
& &-\sum_{i=1}^d\int_{D}\hat{b}_{k,i}(x)\frac{\partial (uv)}{\partial x_i}dx-\int_{D}c_k(x)u(x)v(x)dx,\\
D({\cal E}_k)&=&H^{1,2}_0(D).
\end{eqnarray*}
By (\ref{new334}), following the argument of \cite[Theorem 4.3,
pages 1030-1031]{Zhang3}, we can show that the weak solution to
the Dirichlet boundary value problem
\begin{eqnarray}\label{Kk} \left\{
  \begin{array}{ll}
   L_k u=0 \ \ \ \ {\rm in}\ \ D \\
    u=f \ \ \ \ \ \ \ {\rm on}\ \  \partial D
  \end{array}
\right.
\end{eqnarray}
is given by
\begin{eqnarray*}
u_k(x)&=& E_x^{Q}\left[\exp\left(\int_0^{\tau_D}(c_k-{\rm div}\,\hat b_k)(X_s)ds\right)f(X_{\tau_D})\right] \\
&=& E_x\left[\exp\left(\int_0^{\tau_D}(\tilde
a^{-1}b)^*(X_s)dM_s-\frac{1}{2}\int_0^{\tau_D}b^*\tilde
a^{-1}b(X_s)ds\right.\right. \nonumber\\
& &\ \ \ \ \left.\left. +\int_0^{\tau_D}(c_k-{\rm div}\,\hat
b_k)(X_s)ds\right)f(X_{\tau_D})\right].\nonumber
\end{eqnarray*}

Denote by $v$ the right-hand side of (\ref{15}). We claim that
\begin{equation}\label{claim}
\lim\limits_{k\rightarrow\infty}u_k(x)=v(x),\ \ {\rm for\ q.e.}\ x\in D.
\end{equation} In fact, define
\begin{eqnarray*}
W_k&:=&\exp\left(\int_0^{\tau_D} (c_k-{\rm div}\hat b_k)(X_s)ds\right)\\
&=&\exp\left(\int_0^{\tau_D}c_k(X_s)ds+N_{\tau_D}^{\hat
b_k^H}-\int_0^{\tau_D}\hat b_k^H(X_s)ds\right),\ k\in \mathbb{N},\\
W&:=&\exp\left(\int_0^{\tau_D}c(X_s)ds+N_{\tau_D}^{\hat
b^H}-\int_0^{\tau_D}\hat b^H(X_s)ds\right).
\end{eqnarray*}
By (\ref{addsa}), (\ref{qe11}) and (\ref{qe12}), we get $W_k\rightarrow W$ in
probability under ${Q}_x$ as $k\rightarrow \infty$ for q.e. $x\in
D$. By (\ref{new333}) and Khasminskii's inequality, we obtain that
for $x\in D$,
\begin{eqnarray}\label{Wk}
\sup_{k\in \mathbb{N}}E_x^{Q}[W_k^2]&=& \sup_{k\in
\mathbb{N}}E_x^{Q}
\left[\exp\left(2\int_0^{\tau_D} (c_k-{\rm div}\,\hat b_k)(X_s)ds\right)\right]\nonumber\\
&\le&\sup_{k\in \mathbb{N}}E_x^{Q}
\left[\exp\left(2\int_0^{\tau_D} g_k(X_s)ds\right)\right]\nonumber\\
 &\leq&\frac{1}{1-2\theta}.
\end{eqnarray}
Hence  $\{W_k\}$ is uniformly integrable under ${Q}_x$ for $x\in
D$. Therefore, (\ref{claim}) holds.

Finally, we show  that $v$ is a weak solution  to problem
(\ref{111}). By (\ref{Wk}), we get
\begin{eqnarray}\label{314}
\sup_{k\in \mathbb{N}}\|u_k\|_{L^2}^2&=&\sup_{k\in \mathbb{N}}\int_D\left(E^{{Q}}_x\left[{W_k}f(X_{\tau_D})\right]\right)^2dx \nonumber\\
&<& \frac{\|f\|_{\infty}^2|D|}{1-2\theta},
\end{eqnarray}
where $|D|$ is the Lebesgue measure of $D$. Since $u_k$ is the
weak solution to problem (\ref{Kk}), we have
$\mathcal{E}_k(u_k,\phi)=0$ for any $\phi\in C_0^\infty(D)$. Then,
$\mathcal{E}_k(u_k,\phi)=0,\ \ \forall\phi\in H_0^{1,2}(D)$. Thus,
we have $\mathcal{E}_k(u_k,u_k-u_1)=0$, which implies that
\begin{equation}\label{uku1}
\mathcal{E}_k(u_k,u_k)=\mathcal{E}_k(u_k,u_1). \end{equation}

Note that $|b|^2$, $|\hat b|^2$ and $c$ are in the Kato class. For
any $0<\varepsilon<1$, there exists a constant $A(\varepsilon)>1$
such that for $1\le i\le d$ and $\eta\in H^{1,2}({\mathbb{R}^d})$
(cf. \cite{Kato}),
\begin{equation}\label{kato}
\int_{\mathbb{R}^d}(b_i^2+\hat b_i^2+|c|)\eta^2dx\le
\varepsilon\int_{\mathbb{R}^d}|\nabla
\eta|^2dx+A(\varepsilon)\int_{\mathbb{R}^d}\eta^2dx.
\end{equation}
By (\ref{kato}), we obtain that for $k\in \mathbb{N}$, $1\le i\le
d$ and $\eta\in H^{1,2}({\mathbb{R}^d})$,
\begin{eqnarray}\label{2100}
& &\int_{\mathbb{R}^d}((\hat b_{k,i})^2+|c_k|)\eta^2dx\nonumber\\
&\le& \int_{\mathbb{R}^d}\left\{\int_{\mathbb{R}^d}[\hat b_i^2(x-y)+|c|(x-y)]J_k(y)dy\right\}\eta^2(x)dx\nonumber\\
&\le&\varepsilon\int_{\mathbb{R}^d}|\nabla
\eta|^2dx+A(\varepsilon)\int_{\mathbb{R}^d}\eta^2dx.
\end{eqnarray}
Then, we obtain by (\ref{uku1})-(\ref{2100}) that for $k\in
\mathbb{N}$,
\begin{eqnarray}\label{kato2}
\frac{\lambda}{2}\|\nabla u_k\|_{L^2}^2
&\leq&\frac{1}{2}\sum_{i,j=1}^d\int_{D}a_{ij}(x)\frac{\partial
u_k}{\partial x_i}\frac{\partial u_k}{\partial x_j}dx\nonumber\\
&=&\frac{1}{2}\sum_{i,j=1}^d\int_{D}a_{ij}(x)\frac{\partial
u_k}{\partial x_i}\frac{\partial u_1}{\partial x_j}dx
-\sum_{i=1}^d\int_{D}b_i(x)\frac{\partial u_k}{\partial
x_i}u_1(x)dx\nonumber\\
& &-\sum_{i=1}^d\int_{D}\hat b_{k,i}(x)\frac{\partial
u_k}{\partial x_i}u_1(x)dx
-\sum_{i=1}^d\int_{D}\hat b_{k,i}(x)u_k(x)\frac{\partial u_1}{\partial x_i}dx\nonumber\\
& &-\int_{D}c_k(x)u_k(x)u_1(x)dx+\sum_{i=1}^d\int_{D}b_i(x)\frac{\partial u_k}{\partial x_i}u_k(x)dx\nonumber\\
& &+2\sum_{i=1}^d\int_{D}\hat b_{k,i}(x)\frac{\partial u_k}{\partial x_i}u_k(x)dx+\int_{D}c_k(x)u_k^2(x)dx\nonumber\\
&\leq& \frac{d^2}{2\lambda}\|\nabla u_k\|_{L^2}\|\nabla u_1\|_{L^2}+2dA^{1/2}(\varepsilon)\|\nabla u_k\|_{L^2}\|u_1\|_{H^{1,2}}\nonumber\\
& &+dA^{1/2}(\varepsilon)\|\nabla u_1\|_{L^2}\|u_k\|_{H^{1,2}}+A(\varepsilon)\|u_k\|_{H^{1,2}}\|u_1\|_{H^{1,2}}\nonumber\\
& &+3d\|\nabla u_k\|_{L^2}(\varepsilon\|\nabla
u_k\|^2_{L^2}+A(\varepsilon)\|u_k\|^2_{L^2})^{1/2}\nonumber\\
& & +(\varepsilon\|\nabla
u_k\|^2_{L^2}+A(\varepsilon)\|u_k\|^2_{L^2}).
\end{eqnarray}

Let $\varepsilon$ be much smaller than $\lambda$. Then, we obtain
by (\ref{314}) and (\ref{kato2}) that
$\sup_{k\in\mathbb{N}}\|\nabla u_k\|_{L^2}<\infty$ and thus
$$
\sup_{k\in\mathbb{N}}\|u_k\|_{H^{1,2}}<\infty.
$$
By taking a subsequence if necessary, we may assume that
$u_k\rightarrow v_1$ weakly in $H^{1,2}(D)$ as
$k\rightarrow\infty$ and that its Cesaro mean
$\{u'_k:=\frac{1}{k}\sum_{l=1}^ku_l,\, k\geq1\}\rightarrow v_2$ in
$H^{1,2}(D)$ as $k\rightarrow\infty$. By (\ref{claim}) and \cite[Proposition III.3.5]{MR92}, we obtain that $v_1(x)=v_2(x)=v(x)$ for q.e. $x\in D$ and
\begin{equation}\label{continuous}
v\ {\rm is\ quasi\textrm{-}continuous\ in}\ D.
\end{equation}

Let
$\phi\in C_0^\infty(D)$. Note that for $l\in\mathbb{N}$,
\begin{eqnarray}\label{Ell}
\mathcal{E}_l(u_l,\phi)&=&\frac{1}{2}\sum_{i,j=1}^d\int_{\mathbb{R}^d}a_{ij}(x)\frac{\partial
u_l}{\partial x_i}\frac{\partial \phi}{\partial x_j}dx
-\sum_{i=1}^d\int_{\mathbb{R}^d}b_i(x)\frac{\partial u_l}{\partial x_i}\phi(x)dx\nonumber\\
& &-\sum_{i=1}^d\int_{\mathbb{R}^d}\hat{b}_{l,i}(x)\frac{\partial
(u_l\phi)}{\partial
x_i}dx-\int_{\mathbb{R}^d}c_l(x)u_l(x)\phi(x)dx.
\end{eqnarray}
By (\ref{kato}) and (\ref{2100}), we find that (cf. \cite[Lemma
2.2(iv)]{Zhang3})
\begin{equation}\label{cero1}
\lim_{k\rightarrow\infty}\frac{1}{2}\sum_{i,j=1}^d
\int_{\mathbb{R}^d}a_{ij}(x)\frac{\partial u'_k}{\partial
x_i}\frac{\partial \phi}{\partial
x_j}dx=\frac{1}{2}\sum_{i,j=1}^d\int_{\mathbb{R}^d}a_{ij}(x)\frac{\partial
v}{\partial x_i}\frac{\partial \phi}{\partial x_j}dx,
\end{equation}
\begin{equation}\label{cero2}\lim_{k\rightarrow\infty}\sum_{i=1}^d\int_{\mathbb{R}^d}b_i(x)\frac{\partial u'_k}{\partial x_i}\phi(x)dx
=\sum_{i=1}^d\int_{\mathbb{R}^d}b_i(x)\frac{\partial v}{\partial
x_i}\phi(x)dx,\end{equation}
\begin{equation}\label{cero3}\lim_{k\rightarrow\infty}\frac{1}{k}\sum_{l=1}^k\sum_{i=1}^d\int_{\mathbb{R}^d}\hat{b}_{l,i}(x)\frac{\partial
(u_l\phi)}{\partial x_i}dx
=\sum_{i=1}^d\int_{\mathbb{R}^d}\hat{b}_{i}(x)\frac{\partial
(v\phi)}{\partial x_i}dx,\end{equation} and
\begin{equation}\label{cero4}\lim_{k\rightarrow\infty}\frac{1}{k}\sum_{l=1}^k\int_{\mathbb{R}^d}c_l(x)u_l(x)\phi(x)dx=\int_{\mathbb{R}^d}c(x)v(x)\phi(x)dx.\end{equation}
Therefore, we obtain by (\ref{EEE}) and (\ref{Ell})-(\ref{cero4})
that
$\mathcal{E}(v,\phi)=\lim_{k\rightarrow\infty}\frac{1}{k}\sum_{l=1}^k\mathcal{E}_l(u_l,\phi)=0$.\hfill\fbox

\subsection{Proof of the
continuity of weak solution} It is well-known that any weak
solution to $Lu=0$ in $D$ has a locally
H\"older continuous version (see \cite{Morrey}, cf. also
\cite{M2}). Denote by $v$ the right-hand side of (\ref{15}) and denote by $v^*$ its continuous version in $D$. We will show below that
\begin{equation}\label{v*}
\lim_{x\rightarrow y,x\in D}v^*(x)=f(y),\ \ \forall y\in \partial D.
\end{equation}
First, we prove an important lemma based on
the Dirichlet heat kernel estimates obtained by Aronson.

Suppose
 $d\ge 2$. Let $p_1>d$ and $q_1>1$ satisfy $\frac{1}{p_1}+\frac{1}{q_1}=1$. Then
$q_1=\frac{p_1}{p_1-1}<\frac{d}{d-1}$. We choose $0<\alpha<1$ such
that
\begin{eqnarray}\label{42}
q_1<\frac{d}{d-\alpha}.
\end{eqnarray}
Let $M_1$ be a constant satisfying
\begin{eqnarray}\label{June1}
e^{|x|}\ge M_1|x|^{(d-\alpha+1)/2},\ \ \forall x\in\mathbb{R}^d.
\end{eqnarray}
 Let $p_2> d/2$ and $q_2>1$
satisfy $\frac{1}{p_2}+\frac{1}{q_2}=1$. Then
$q_2=\frac{p_2}{p_2-1}<\frac{d}{d-2}$. We choose $\beta$
satisfying
\begin{eqnarray}\label{new42}\frac{d}{2}-1<\beta<\frac{d}{2q_2}.\end{eqnarray}
Let $M_2$ be a constant satisfying
\begin{eqnarray}\label{June11}
e^{|x|}\ge M_2|x|^{\beta},\ \ \forall x\in\mathbb{R}^d,
\end{eqnarray}
and let $M_3$ be a constant satisfying
$$
e^{|x|}\ge M_3|x|^{5/8},\ \ \forall x\in\mathbb{R}^d.
$$
We denote by $\varsigma$ the diameter of $D$ as above. By
(\ref{42}) and (\ref{new42}), we find that
$$
\int_0^\varsigma r^{d-q_1(d-\alpha)-1}dr<\infty\ \ {\rm and}\ \
\int_0^\varsigma r^{d-2\beta q_2-1}dr<\infty.
$$

Denote
$$h(t,x,y)=\frac{\sigma_1}{t^{d/2}}e^{-\frac{\sigma_2|x-y|^2}{t}},\ \ (t,x,y)\in (0,\infty)\times D\times D.$$
Then, we obtain by (\ref{213}) that
$$p(t,x,y)\leq h(t,x,y),\ \ (t,x,y)\in (0,\infty)\times D\times D.
$$

\begin{lem}\label{lemma23} Let $\mu$ be a vector field on $\mathbb{R}^d$ and $\nu$ be a function on $\mathbb{R}^d$ such that $\mu,\nu\in
C^{\infty}(\mathbb{R}^d)$.

(i) Suppose $d\ge 2$, $p_1>d$ and $p_2> d/2$. Then, for $t>0$ and $x\in D$,
\begin{eqnarray*}&&\left|\int_{y\in D}h(t,x,y){\rm div}\,\mu(y)dy\right|\nonumber\\
&\le
&\frac{2\sigma_1}{\sigma_2^{(d-\alpha-1)/2}M_1t^{(1+\alpha)/2}}\left(\int_0^\varsigma
r^{d-q_1(d-\alpha)-1}dr\right)^{1/q_1}\left(\int_{y\in
D}|\mu(y)|^{p_1}dy\right)^{1/p_1}
\end{eqnarray*}
and
\begin{eqnarray*}&&\left|\int_{y\in D}h(t,x,y)\nu(y)dy\right|\nonumber\\
&\le &\frac{\sigma_1}{\sigma_2^\beta M_2t^{d/2-\beta}}\left(\int_0^\varsigma
r^{d-2\beta q_2-1}dr\right)^{1/q_2}\left(\int_{y\in
D}|\nu(y)|^{p_2}dy\right)^{1/p_2}.
\end{eqnarray*}

(ii)  Suppose $d=1$. Then, for $t>0$ and $x\in D$,
\begin{eqnarray*}&&\left|\int_{y\in D}h(t,x,y){\rm div}\,\mu(y)dy\right|\le \frac{2^{3/2}\sigma_1\sigma_2^{3/8}\varsigma^{1/4}}{M_3t^{7/8}}\left(\int_{y\in
D}|\mu(y)|^2dy\right)^{1/2}
\end{eqnarray*}
and
\begin{eqnarray*}\left|\int_{y\in D}h(t,x,y)\nu(y)dy\right|\le \frac{\sigma_1}{t^{1/2}}\int_{y\in D}|\nu(y)|dy.
\end{eqnarray*}
\end{lem}

\noindent {\bf Proof.} We only prove (i). The proof of (ii) is
similar so we omit it here.

By (\ref{June1}), we get
\begin{eqnarray*}
& &\left|\int_{y\in D}h(t,x,y){\rm div}\,\mu(y)dy\right|\nonumber\\
&=&\left|\int_{y\in D}\langle\nabla_y h(t,x,y),\mu(y)\rangle dy\right|\nonumber\\
&\le&\int_{y\in
D}\frac{\sigma_1}{t^{d/2}e^{\sigma_2|x-y|^2/t}}\frac{2\sigma_2|x-y|}{t}|\mu
(y)|dy\\
&\le&\frac{2\sigma_1\sigma_2}{t^{d/2+1}}\int_{y\in D}\frac{|x-y|}{M_1(\sigma_2|x-y|^2/t)^{(d-\alpha+1)/2}}|\mu(y)|dy\nonumber\\
&\le&\frac{2\sigma_1}{\sigma_2^{(d-\alpha-1)/2}M_1t^{(1+\alpha)/2}}\int_{y\in D}\frac{|\mu(y)|}{|x-y|^{d-\alpha}}dy\nonumber\\
&\le&\frac{2\sigma_1}{\sigma_2^{(d-\alpha-1)/2}M_1t^{(1+\alpha)/2}}\left(\int_{y\in
D}\frac{1}{|x-y|^{q_1(d-\alpha)}}dy\right)^{1/q_1}\left(\int_{y\in
D}|\mu(y)|^{p_1}dy\right)^{1/p_1}
\nonumber\\
&\le&\frac{2\sigma_1}{\sigma_2^{(d-\alpha-1)/2}M_1t^{(1+\alpha)/2}}\left(\int_0^\varsigma
r^{d-q_1(d-\alpha)-1}dr\right)^{1/q_1}\left(\int_{y\in
D}|\mu(y)|^{p_1}dy\right)^{1/p_1}.
\end{eqnarray*}

By (\ref{June11}), we get
\begin{eqnarray*}
& &\left|\int_{y\in D}h(t,x,y)\nu(y)dy\right|\nonumber\\
&\le&\int_{y\in D}\frac{\sigma_1}{t^{d/2}e^{\sigma_2|x-y|^2/t}}|\nu
(y)|dy\\
&\le&\int_{y\in D}\frac{\sigma_1}{M_2t^{d/2}(\sigma_2|x-y|^2/t)^{\beta}}|\nu(y)|dy\nonumber\\
&=&\frac{\sigma_1}{\sigma_2^\beta M_2t^{d/2-\beta}}\int_{y\in D}\frac{|\nu(y)|}{|x-y|^{2\beta}}dy\nonumber\\
&\le&\frac{\sigma_1}{\sigma_2^\beta M_2t^{d/2-\beta}}\left(\int_{y\in
D}\frac{1}{|x-y|^{2\beta q_2}}dy\right)^{1/q_2}\left(\int_{y\in
D}|\nu(y)|^{p_2} dy\right)^{1/p_2}
\nonumber\\
&\le&\frac{\sigma_1}{\sigma_2^\beta M_2t^{d/2-\beta}}\left(\int_0^\varsigma
r^{d-2\beta q_2-1}dr\right)^{1/q_2}\left(\int_{y\in
D}|\nu(y)|^{p_2}dy\right)^{1/p_2}.
\end{eqnarray*}\hfill\fbox

\begin{rem}
In \cite{Cho}, Cho, Kim and Park established very nice sharp
two-sided estimates on Dirichlet heat kernels. Under the
additional assumption that $D$ is a $C^{1,\alpha}$-domain
($0<\alpha\le1$) satisfying the connected line condition and each
$a_{ij}$, $1\le i,j\le d$, is Dini continuous, by \cite[Theorem
1.1]{Cho},  for each $T>0$, there exist positive constants $c_1$
and $c_2$ such that for $(t,x,y)\in(0,T)\times D\times D$,
\begin{equation}\label{gauss} p(t,x,y)\le
\left(1\wedge\frac{\rho(x)}{\sqrt{t}}\right)\left(1\wedge\frac{\rho(y)}{\sqrt{t}}\right)
\frac{c_1}{t^{d/2}}e^{-\frac{c_2|x-y|^2}{t}}
\end{equation}
and
\begin{equation}\label{gauss2} |\nabla_yp(t,x,y)|\le
\left(1\wedge\frac{\rho(x)}{\sqrt{t}}\right)\frac{c_1}{t^{(d+1)/2}}e^{-\frac{c_2|x-y|^2}{t}},
\end{equation}
where $\rho(x):={dist}(x,\partial D)$.

By virtue of (\ref{gauss}) and (\ref{gauss2}), we can obtain
estimates for $p(t,x,y)$ similar to those for $h(t,x,y)$ given as
in Lemma \ref{lemma23}. These estimates for $p(t,x,y)$ or
$h(t,x,y)$ make it possible to handle the case when Meyers's
$L^p$-estimate is not available.
\end{rem} \vskip 0.5cm \noindent {\bf Proof of the continuity of
weak solution at the boundary.} \vskip 0.3cm

By (\ref{continuous}), we have $v^*(x)=v(x)$ for q.e. $x\in D$. Note that for $x\in D$,
\begin{eqnarray*}
v(x)&=&E^Q_x\left[\exp\left(\int_0^{\tau_D} c(X_s)ds+N^{{\hat
b}^H}_{{\tau_D}}-\int_0^{\tau_D} {\hat
b}^H(X_s)ds\right)f(X_{{\tau_D}})\right]\\
&=&E_x^{Q}[f(X_{\tau_{D}})]+E_x^Q[f(X_{\tau_{D}})(e^{A_{\tau_D}}-1)],
\end{eqnarray*}
where $A_{t}:=\int_0^{t}c(X_s)ds+N_{t}^{\hat b^H}-\int_0^{t}\hat
b^H(X_s)ds$, $t\ge 0$. By Lemma \ref{05}, to prove (\ref{v*}), it
suffices to show that there exists an exceptional set $F\subset D$
such that
\begin{equation}\label{40}
\lim_{x\rightarrow y,x\in D\backslash
F}E_x^{Q}[f(X_{\tau_{D}})(e^{A_{\tau_D}}-1)]=0,\ \ \forall
y\in\partial D.
\end{equation}

For $t>0$ and $x\in D$, we have
\begin{eqnarray*}E_x^{Q}[f(X_{\tau_{D}})(e^{A_{\tau_D}}-1)]
&=&E_x^{Q}[f(X_{\tau_{D}})(e^{A_{\tau_D}}-1);\tau_D\leq
t]\\
&
&+E_x^{Q}[f(X_{\tau_{D}})(e^{A_{\tau_D}}-1);\tau_D>t].\end{eqnarray*}
By (\ref{ZtauD}), there exists an exceptional set $F_1\subset D$
such that
$$
\sup_{x\in D\backslash F_1}E_x^{Q}[\exp({A_{\tau_D}})]=\sup_{x\in
D\backslash F_1}E_x[Z_{\tau _D}]\le\frac{1}{1-\theta}.
$$
Then, we obtain by the strong Markov property that for q.e. $x\in D$,
\begin{eqnarray}\label{great}
& &\left|E_x^{Q}[f(X_{\tau_{D}})(e^{A_{\tau_D}}-1);\tau_D>t]\right|\nonumber\\
&\leq&\|f\|_\infty\left\{{Q}_x(\tau_D>t)+ E_x^{Q}[e^{A_{\tau_D}};\tau_D>t]\right\}\nonumber\\
&\le& \|f\|_\infty
\left\{{Q}_x(\tau_D>t)+\frac{E_x^{Q}[e^{A_t};\tau_D>t]}{1-\theta}\right\},\ \ \forall t>0.
\end{eqnarray}
By Lemma \ref{05}, following the argument of
\cite[(2.28)]{Chen}, we get
\begin{eqnarray}\label{great21}\lim_{x\rightarrow y,x\in D}{Q}_x(\tau_D>t)=0,\ \ \forall t>0,\forall y\in\partial D.\end{eqnarray}
By   (\ref{addsa}), (\ref{qe11}),  (\ref{qe12}),  (\ref{Wk}) and
Fatou's lemma,  there exists an exceptional set $F_2\subset D$
such that  for every $t>0$,
\begin{eqnarray}\label{great33}\sup_{x\in D\backslash F_2}E_x^{Q}[e^{2A_t};\tau_D>t]\leq\sup_{x\in D\backslash F_2}\sup_{k\in \mathbb{N}}E_x^{Q}\left[e^{2\int_0^{\tau_D}g_k(X_s)ds}ds\right]
\le\frac{1}{1-2\theta}.\end{eqnarray} Thus, we obtain by
(\ref{great})-(\ref{great33}) that there exists an exceptional set
$F_3\subset D$ satisfying
$$\lim_{x\rightarrow y,x\in D\backslash F_3}E_x^{Q}[f(X_{\tau_{D}})(e^{A_{\tau_D}}-1);\tau_D>t]=0,\ \ \forall t>0,\forall y\in\partial D.$$
Therefore, to prove (\ref{40}), it suffices to show that  there
exists an exceptional set $F_4\subset D$ such that
\begin{equation}\label{xcv}
\lim_{t\downarrow0}\sup_{x\in D\backslash
F_4}E_x^{Q}[f(X_{\tau_{D}})(e^{A_{\tau_D}}-1);\tau_D\le t]=0.
\end{equation}

By  (\ref{addsa}), (\ref{qe11}),  (\ref{qe12}) and Fatou's lemma,
there exists an exceptional set $F_4\subset D$ such that  for
every $t>0$,
\begin{eqnarray*}
&&\sup_{x\in D\backslash F_4}|E_x^{Q}[f(X_{\tau_{D}})(e^{A_{\tau_D}}-1);\tau_D\leq t]|\\
&\leq&\|f\|_\infty\sup_{x\in D\backslash
F_4}\liminf_{k\rightarrow\infty}E_x^{Q}
\left[\left|e^{\int_0^{\tau_D}(c_k-{\rm div}\,\hat b_k)(X_s)ds}-1\right|;\tau_D\leq t\right]\\
&\leq&\|f\|_\infty\left\{\sup_{x\in D}\limsup_{k\rightarrow\infty}E_x^{Q}\left[e^{\int_0^{\tau_D}g_k(X_s)ds}-1;\tau_D\leq t\right]\right.\\
& & \left.+\sup_{x\in
D}\limsup_{k\rightarrow\infty}E_x^{Q}\left[\left(1-e^{\int_0^{\tau_D}(c_k-{\rm
div}\,\hat b_k-g_k)(X_s)ds}\right);\tau_D\leq t\right]\right\}\\
&\leq&\|f\|_\infty\left\{\sup_{x\in
D}\limsup_{k\rightarrow\infty}E_x^{Q}\left[e^{\int_0^{t\wedge\tau_D}g_k(X_s)ds}-1
\right]\right.\\
& & \left.+\sup_{x\in
D}\limsup_{k\rightarrow\infty}E_x^{Q}\left[\left(1-e^{\int_0^{t\wedge\tau_D}(c_k-{\rm
div}\,\hat b_k-g_k)(X_s)ds}\right)\right]\right\}.
\end{eqnarray*}
By Lemma \ref{lemma23} and Khasminskii's inequality, we get
$$\lim_{t\downarrow0}\sup_{x\in D}\sup_{k\in\mathbb{N}}E_x^{Q}\left[e^{\int_0^{t\wedge\tau_D} g_k(X_s)ds}\right]= 1.$$
Hence, to prove (\ref{xcv}), we need only show that
$$
\lim_{t\downarrow0}\inf_{x\in
D}\inf_{k\in\mathbb{N}}E_x^{Q}\left[e^{\int_0^{t\wedge\tau_D}(c_k-{\rm
div}\,\hat b_k-g_k)(X_s)ds}\right]\ge1.
$$
Further, by Jensen's inequality, we need only show that
$$
\lim_{t\downarrow0}\sup_{x\in
D}\sup_{k\in\mathbb{N}}E_x^{Q}\left[{\int_0^{t\wedge\tau_D}(g_k-c_k+{\rm
div}\,\hat b_k)(X_s)ds}\right]=0.
$$

By Lemma \ref{lemma23}, we obtain that
\begin{eqnarray*}
& &\sup_{x\in
D}\sup_{k\in\mathbb{N}}E_x^{Q}\left[{\int_0^{t\wedge\tau_D}(g_k-c_k+{\rm
div}\,\hat b_k)(X_s)ds}\right]\\
&=&\sup_{x\in D}\sup_{k\in\mathbb{N}}\int_0^t\int_{y\in
D}p(s,x,y)(g_k-c_k+{\rm div}\,\hat b_k)(y)dyds\\
&\leq&\sup_{x\in D}\sup_{k\in\mathbb{N}}\int_0^t\int_{y\in
D}h(s,x,y)(g_k-c_k+{\rm div}\,\hat b_k)(y)dyds\\
&\rightarrow&0\ {\rm as}\ t\downarrow0.
\end{eqnarray*}
The proof is complete.\hfill\fbox

\subsection{Proof of the uniqueness of continuous weak solutions}

In this subsection, we will prove that there exists a unique
continuous weak solution to problem (\ref{111}).

Let $u_1$ be a weak solution to problem  (\ref{111}) such that
$u_1$ is continuous on $\overline{D}$. We have Fukushima's
decomposition
\begin{eqnarray}\label{july11}
u_1(X_t)- u_1(X_0)&=&M_t^{u_1}+N_t^{u_1}\nonumber\\
&=&\int_0^t\nabla u_1(X_s)dM_s+N_t^{u_1},\ \
t<\tau_D.\end{eqnarray} We claim that
\begin{eqnarray}\label{309}
N_t^{u_1}&=& -\sum_{i=1}^{d}\int_0^tb_i(X_s)\frac{\partial u_1}{\partial x_i}(X_s)ds-\int_0^tu_1(X_s)c(X_s)ds \nonumber\\
& &-\int_0^t u_1(X_s)dN_s^{\hat b^H}+\int_0^tu_1(X_s)\hat
b^H(X_s)ds,\nonumber\\
& &\ \ \ \ \ \ \ \ \ \ \ \ \ \ \ \ t<\tau_D,\ P_x-{\rm a.s.}\ {\rm for\ q.e.}\ x\in D,
\end{eqnarray}
where the third term of (\ref{309}) is a Nakao integral (we refer
the readers to \cite[Definition 2.4]{CMS} and
\cite[Definition 3.1]{N} for the definition).

Let $\{D_n\}$ be a sequence of increasing open subsets of $
\mathbb{R}^d$ satisfying $D=\cup_{n\in\mathbb{N}}D_n$ and
$\overline{D_n}\subset D_{n+1}$ for each $n$. We choose a sequence
$\{u^{(n)}\subset H^{1,2}_0(D)\cap {\cal B}_b(D_n)\}$ satisfying
$u_1=u^{(n)}$ on $D_n$ for each $n$. To prove (\ref{309}) it
suffices to show that for any $n\in\mathbb{N}$,
\begin{eqnarray}\label{nna}
N_t^{u^{(n)}}&=& -\sum_{i=1}^{d}\int_0^tb_i(X_s)\frac{\partial u^{(n)}}{\partial x_i}(X_s)ds-\int_0^tu^{(n)}(X_s)c(X_s)ds \nonumber\\
& &-\int_0^t u^{(n)}(X_s)dN_s^{\hat b^H}+\int_0^tu^{(n)}(X_s)\hat
b^H(X_s)ds, \nonumber\\
& &\ \ \ \ \ \ \ \ \ \ \ \ \ \ \ \ t<\tau_{D_n},\ P_x-{\rm a.s.}\ {\rm for\ q.e.}\ x\in D.
\end{eqnarray}
Denote by $C^{(n)}_t$ the right hand side of (\ref{nna}). By \cite[Theorem 5.2.7]{oshima}, following the argument of the proof of \cite[Theorem 2.2]{N}, we find that to prove (\ref{nna}) it suffices to show that for each $n$,
\begin{eqnarray}\label{CVB}
\lim_{t\downarrow0}\frac{1}{t}E_{\phi\cdot
dx}[N_t^{u^{(n)}}]=\lim_{t\downarrow0}\frac{1}{t}E_{\phi\cdot
dx}[C^{(n)}_t],\ \ \forall \phi\in H^{1,2}_0(D_n)\cap {\cal
B}_b(D_n).
\end{eqnarray}

We fix an $n\in \mathbb{N}$ and  $\phi\in H^{1,2}_0(D_n)\cap {\cal
B}_b(D_n)$. By (\ref{EEE}), (\ref{EEF}) and (\ref{2002}), we get
\begin{eqnarray}\label{317}
\mathcal{E}^0(u^{(n)},\phi)&=&\mathcal{E}(u^{(n)},\phi)+\sum_{i=1}^d\int_{D}b_i(x)\frac{\partial u^{(n)}}{\partial x_i}\phi(x)dx\nonumber\\
& &+\sum_{i=1}^d\int_{D}\hat{b}_i(x)\frac{\partial (u^{(n)}\phi)}{\partial x_i}+\int_{D}c(x)u^{(n)}(x)\phi(x)dx,\nonumber\\
&=&\sum_{i=1}^{d}\int_Db_i(x)\frac{\partial u^{(n)}}{\partial
x_i}\phi(x)dx+\int_Dc(x)u^{(n)}(x)\phi(x)dx\nonumber\\
& &-{\cal E}_1^0(\hat b^H,u^{(n)}\phi).
\end{eqnarray}
We have
\begin{eqnarray}\label{3001}
\lim_{t\downarrow0}\frac{1}{t}E_{\phi\cdot dx}[N_t^{u^{(n)}}]&=&\lim_{t\downarrow0}\frac{1}{t}E_{\phi\cdot dx}[u^{(n)}(X_t)-u^{(n)}(X_0)-M_t^{u^{(n)}}]\nonumber\\
&=&\lim_{t\downarrow0}\frac{1}{t}\int_{D}E_x[u^{(n)}(X_t)-u^{(n)}(X_0)]\phi(x)dx\nonumber\\
&=&-{\cal E}^0(u^{(n)},\phi)
\end{eqnarray}
and
\begin{eqnarray}\label{3002}
& &\lim_{t\downarrow0}\frac{1}{t}E_{\phi\cdot
dx}\left[-\sum_{i=1}^{d}\int_0^tb_i(X_s)\frac{\partial
u^{(n)}}{\partial x_i}(X_s)ds
-\int_0^tu^{(n)}(X_s)c(X_s)ds\right.\nonumber\\
& &\left. \ \ \ \ \ \ \ \ +\int_0^tu^{(n)}(X_s)\hat
b^H(X_s)ds\right]\nonumber\\
&=&-\sum_{i=1}^{d}\int_Db_i(x)\frac{\partial u^{(n)}}{\partial
x_i}\phi(x)dx-\int_{D}c(x)u^{(n)}(x)\phi(x)dx\nonumber\\
& &+\int_{D}\hat b^H(x)u^{(n)}(x)\phi(x)dx.
\end{eqnarray}
By \cite[Remark 2.5]{CMS}, we get
\begin{equation}\label{3005}
\lim_{t\downarrow0}\frac{1}{t}E_{\phi\cdot
dx}\left[-\int_0^tu^{(n)}(X_s)dN_s^{\hat b^H}\right]={\cal
E}^0(\hat b^H, u^{(n)}\phi).
\end{equation}
Then, (\ref{CVB}) holds by (\ref{317})-(\ref{3005}). Thus,
(\ref{nna}) and hence (\ref{309}) hold.

By (\ref{july11}) and (\ref{309}), we obtain that
\begin{eqnarray}\label{beiyong}
& &u_1(X_t)-u_1(X_0) \nonumber\\ &=& \int_0^t\nabla u_1(X_s)dM_s-\sum_{i=1}^{d}\int_0^tb_i(X_s)\frac{\partial u_1}{\partial x_i}(X_s)ds \nonumber\\
& &-\int_0^tu_1(X_s)c(X_s)ds-\int_0^tu_1(X_s)dN_s^{\hat
b^H}+\int_0^tu_1(X_s)\hat b^H(X_s)ds, \nonumber\\
& &\ \ \ \ \ \ \ \ \ \ \ \ \ \ \ \ t<\tau_{D},\ P_x-{\rm a.s.}\ {\rm for\ q.e.}\ x\in D.
\end{eqnarray}
We now prove that for $t<\tau_D$,
\begin{equation}\label{318}
d(u_1(X_t)Z_t)=u_1(X_t)Z_t(\tilde a^{-1}b)^*(X_t)dM_t+Z_t\nabla
u_1(X_t)dM_t,
\end{equation}
$P_x-{\rm a.s.}\ {\rm for\ q.e.}\ x\in D$, where $Z_t$ is defined as in (\ref{311}).

For $k\in \mathbb{N}$ and $t>0$, we define
\begin{eqnarray}\label{bew}
V^k_t&:=& \int_0^t\nabla u_1(X_s)dM_s-\sum_{i=1}^{d}\int_0^tb_i(X_s)\frac{\partial u_1}{\partial x_i}(X_s)ds\nonumber\\
& &-\int_0^tu_1(c_k-{\rm div}\,\hat b_k)(X_s)ds
\end{eqnarray}
and
\begin{eqnarray*}
Z^k_t&:=&\exp\left(\int_0^t(\tilde
a^{-1}b)^*(X_s)dM_s-\frac{1}{2}\int_0^{t}
b^*\tilde a^{-1}b(X_s)ds\right.\\
& &\ \ \left. +\int_0^t(c_k-{\rm div}\,\hat b_k)(X_s)ds\right).
\end{eqnarray*}
Then,
$$dZ^k_t=Z^k_t(\tilde
a^{-1}b)^*(X_t)dM_t+Z^k_t(c_k-{\rm div}\,\hat b_k)(X_t)dt.$$ Note
that both $\{V^k_t\}$ and $\{Z^k_t\}$ are semi-martingales.
Applying { Ito's formula}, we obtain that
\begin{eqnarray*}
d(V^k_tZ^k_t)&=&V^k_tZ^k_t(\tilde
a^{-1}b)^*(X_t)dM_t+Z^k_t\nabla u_1(X_t)dM_t\nonumber\\
& &+Z_t^k(V^k_t-u_1(X_t))(c_k-{\rm div}\,\hat b_k)(X_t)dt.
\end{eqnarray*}
Further, applying {Ito's formula} to $Z^k_t$, we get
\begin{eqnarray}\label{ito2}
& &d((V^k_t+u_1(X_0))Z^k_t)\nonumber\\
&=&V^k_tZ^k_t(\tilde
a^{-1}b)^*(X_t)dM_t+Z^k_t\nabla u_1(X_t)dM_t\nonumber\\
& &+Z_t^k(V^k_t-u_1(X_t))(c_k-{\rm div}\,\hat b_k)(X_t)dt\nonumber\\
& &+u_1(X_0)Z^k_t(\tilde
a^{-1}b)^*(X_t)dM_t+u_1(X_0)Z^k_t(c_k-{\rm div}\,\hat b_k)(X_t)dt\nonumber\\
&=&(V^k_t+u_1(X_0))Z^k_t(\tilde
a^{-1}b)^*(X_t)dM_t+Z^k_t\nabla u_1(X_t)dM_t\nonumber\\
& &+Z_t^k(V^k_t-(u_1(X_t)-u_1(X_0)))(c_k-{\rm div}\,\hat
b_k)(X_t)dt.
\end{eqnarray}
By  (\ref{addsa}), (\ref{qe11}), (\ref{beiyong}), (\ref{bew})  and \cite[Theorem 2.7]{CMS}, there exists a subsequence $\{k_l\}$ such that  $V^{k_l}_t\rightarrow u_1(X_t)-u_1(X_0)$, $t<\tau_D$, $P_x-{\rm a.s.}\ {\rm for\ q.e.}\ x\in D$ as $l\rightarrow\infty$. Therefore, (\ref{318}) holds by (\ref{ito2}).

By (\ref{318}), we know that $\{u_1(X_{t\wedge
\tau_D})Z_{t\wedge\tau_D},t\geq0\}$ is a ${P}_x$-local martingale
for q.e. $x\in D$. We claim that $\{Z_{t\wedge \tau_D},t\geq0\}$
is ${P}_x$-uniformly integrable for q.e. $x\in D$. Write
$$Z_{t\wedge \tau_D}=Z_{\tau_D}{1}_{\{\tau_D\le
t\}}+Z_t{1}_{\{\tau_D> t\}}.$$ By (\ref{ZtauD}),
$\{Z_{\tau_D}{1}_{\{\tau_D\le t\}}, t\ge 0\}$ is ${P}_x$-uniformly
integrable for q.e. $x\in D$. We now show that $\{Z_t{1}_{\{\tau_D> t\}}, t\ge 0\}$
is ${P}_x$-uniformly integrable for q.e. $x\in D$. Note that  for q.e. $x\in D$,
\begin{eqnarray*}
Z_t{1}_{\{\tau_D> t\}}&\leq&{1}_{\{\tau_D>
t\}}\exp\left(\int_0^{\tau_D}(\tilde
a^{-1}b)^*(X_s)dM_s-\frac{1}{2}\int_0^{\tau_D}
b^*\tilde a^{-1}b(X_s)ds\right.\\
& &\ \ \left. +\int_0^{\tau_D} g(X_s)ds\right)\\
&:=&{1}_{\{\tau_D> t\}}Z^g_{\tau_D}.
\end{eqnarray*}
Hence it suffices to show that $\{{1}_{\{\tau_D>
t\}}Z^g_{\tau_D},t\ge 0\}$ is ${P}_x$-uniformly integrable for $x\in D$.

By the strong Markov property, we get
\begin{eqnarray}\label{3119}
{1}_{\{\tau_D> t\}}E_x[Z^g_{\tau_D}|\mathcal{F}_t]&=&{1}_{\{\tau_D> t\}}Z^g_tE_{X_t}[Z^g_{\tau_D}]\nonumber\\
&\geq&{1}_{\{\tau_D> t\}}Z^g_t\inf_{x\in D} E_{x}[Z^g_{\tau_D}]\nonumber\\
&=&{1}_{\{\tau_D> t\}}Z^g_t\inf_{x\in D} E^{{Q}}_{x}\left[\exp\left(\int_0^{\tau_D}g(X_s)ds\right)\right]\nonumber\\
&\geq&{1}_{\{\tau_D> t\}}Z^g_t.
\end{eqnarray}
By (\ref{old334}) and (\ref{3119}), we obtain that
$\{{1}_{\{\tau_D> t\}}Z^g_{\tau_D},t\ge 0\}$ is ${P}_x$-uniformly
integrable for $x\in D$. Therefore $\{Z_{t\wedge\tau_D},t\geq0\}$ is
${P}_x$-uniformly integrable for q.e. $x\in D$. Since $u_1$ is
bounded continuous, we find that $\{u_1(X_{t\wedge
\tau_D})Z_{t\wedge\tau_D},t\geq0\}$ is a ${P}_x$-martingale for
q.e. $x\in D$. Thus,
$$u_1(x)=E_x[u_1(X_{t\wedge \tau_D})Z_{t\wedge\tau_D}],\ \ {\rm for\ q.e.}\ x\in D.$$
Letting $t\rightarrow\infty$, we obtain that
$$u_1(x)=E_x[f(X_{\tau_D})Z_{\tau_D}],\ \ {\rm for\ q.e.}\ x\in D,$$
which proves the uniqueness. \hfill\fbox

\section{Probabilistic Representation of Non-symmetric Semigroup}\setcounter{equation}{0}

In this section, we will use some techniques of Section 2 to give
a probabilistic representation of the non-symmetric semigroup
$\{T_t\}_{t\ge 0}$ associated with the operator $L$ defined by
(\ref{12}). The obtained result (see Theorem \ref{main} below)
generalizes \cite[Theorem 3.4]{Zhang1}, which is the first result
on the probabilistic representation of semigroups with $\hat
b\not=0$, from the case of symmetric diffusion matrix $A$ to the
non-symmetric case. The methods and techniques of this paper can
be applied also to some other problems such as the mixed boundary
value problem, Dirichlet problem of semilinear elliptic PDEs with
singular coefficients, etc. (cf. \cite{Zhang5,{Zhang4}}). We will
consider them in future work.

Throughout this section, we let $D$ be an open subset of
$\mathbb{R}^d$, which need not be bounded. Suppose that
$A(x)=(a_{ij}(x))_{i,j=1}^d$ is a Borel measurable matrix-valued
function  on $D$ satisfying (\ref{13}) and (\ref{14});
$b=(b_1,\dots,b_d)^*$ and $\hat{b}=(\hat{b}_1,\dots,\hat{b}_d)^*$
are Borel measurable $\mathbb{R}^d$-valued functions on $D$ and
$c$ is a Borel measurable function on $D$ satisfying $|b|^2\in
L^{p\vee 1}(D;dx), |\hat b|^2\in L^{p\vee 1}(D;dx)$ and $c\in
L^{p\vee 1}(D;dx)$ for some constant $p>d/2$. Let $L$ and $({\cal
E},D(\cal{E}))$ be defined as in (\ref{12}) and (\ref{EEE}),
respectively. Since $|b|^2$, $|\hat b|^2$ and $c$ are in the Kato
class, there exists a constant $\gamma>0$ such that $({\cal
E}_{\gamma},D(\cal{E}))$ is a coercive closed form on $L^2(D;dx)$
(cf. \cite[page 329]{Zhang1}). Hence there exits a (unique)
strongly continuous semigroup $\{T_t\}_{t\ge 0}$ on $L^2(D;dx)$
which is associated with $({\cal E},D(\cal{E}))$. Denote by
$({\cal L}, D({\cal L}))$  the generator of $\{T_t\}_{t\ge 0}$.
 Clearly ${\cal L}$ is formally given by $L$. Denote by
$\{\hat{T}_t\}_{t\ge 0}$ the dual semigroup of $\{T_t\}_{t\ge 0}$
on $L^2(D;dx)$.

We define the Dirichlet form $({\cal E}^0,D({\cal E}^0))$ as in
(\ref{EEF}). Let $X=((X_t)_{t\ge 0},(P_x)_{x\in \mathbb{R}^d})$
and $\hat X=(X_t)_{t\ge 0},(\hat P_x)_{x\in \mathbb{R}^d})$ be the
Markov process and dual Markov process associated with the
Dirichlet form $({\cal E}^0,D({\cal E}^0))$ given by (\ref{EEF}),
respectively. Let $M_t$, $(\tilde a_{ij})_{i,j=1}^d$, $v^H$, etc.
be defined the same as in Section 1. Denote by $m$ the Lebesgue
measure $dx$ on $\mathbb{R}^d$. Now we can state the main result
of this section.

\begin{thm}\label{main} For any $f,g \in L^2(D;dx)$, we have
\begin{eqnarray} \label{28}
& &\int_{D} f(x)T_tg(x)dx\nonumber\\
&=&E_m\left[f(X_0)g(X_t)\exp\left(\int_0^t(\tilde a^{-1}b)^*(X_s)dM_s-\frac{1}{2}\int_0^tb^*\tilde{a}^{-1}b(X_s)ds\right.\right.\nonumber\\
& &\ \ \ \ \left.\left. +\int_0^t c(X_s)ds+N^{{\hat
b}^H}_t-\int_0^t {\hat b}^H(X_s)ds\right); t<\tau_D\right].
\end{eqnarray}
\end{thm}

\noindent {\bf Proof.} By (\ref{30}), similar to \cite[Theorem 2.1]{Zhang1}, we can prove
the following lemma on integrability of functionals of Dirichlet
processes.
\begin{lem}\label{432}
Suppose $f\in L^{r\vee 1}(D;dx)$ for some $r>d/2$ and $T>0$. Then,
there exists a constant $\varrho_1>0$ depending on $f$, $r$ and
$T$ such that for any $0\le t\le T$,
$$
\sup_{x\in D}E_x\left[\exp\left(\int_0^tf(X_s)ds\right);
t<\tau_D\right]\le \varrho_1e^{\varrho_1t},
$$
and
$$
\sup_{x\in D}\hat{E}_x\left[\exp\left(\int_0^tf(X_s)ds\right);
t<\tau_D\right]\le \varrho_1e^{\varrho_1t}.
$$

\end{lem}

We divide the proof of Theorem \ref{main} into three cases.

\noindent {\bf Case 1: $\hat b=0$}.

For $g\in {\cal B}_b(D)$, we define
\begin{eqnarray*}
& &P_tg(x):=E_x\left[\exp\left(\int_0^t(\tilde
a^{-1}b)^*(X_s)dM_s\right.\right.\nonumber\\
& &\ \ \ \ \ \ \ \left.\left.-\frac{1}{2}\int_0^t b^*
\tilde{a}^{-1} b(X_s)ds + \int_0^t c(X_s)ds \right)g(X_t);
t<\tau_D\right].
\end{eqnarray*}
Clearly $\{P_t\}_{t\ge 0}$ is a well-defined semigroup. We now
show that $\{P_t\}_{t\ge 0}$ extends to a strongly continuous
semigroup on $L^2(D;dx)$, which will be also denoted by
$\{P_t\}_{t\ge 0}$.

In fact, for any $g\in L^2(D;dx)$, we obtain by Lemma \ref{432}
that
 \begin{eqnarray}\label{223}
 &&\int_{D}(P_tg(x))^2dx\nonumber\\
 &=&\int_{D}\left(E_x\left[\exp\left(\int_0^t(\tilde a^{-1}b)^*(X_s)dM_s-\int_0^t b^* \tilde{a}^{-1} b(X_s)ds\right) \right.\right. \nonumber \\
 & &\ \ \ \ \left.\left. \cdot\exp\left(\frac{1}{2}\int_0^t b^* \tilde{a}^{-1} b(X_s)ds + \int_0^t c(X_s)ds \right)g(X_t);
t<\tau_D\right]\right)^2dx \nonumber \\
 &\leq& \int_{D}E_x\left[\exp\left(2\int_0^t(\tilde a^{-1}b)^*(X_s)dM_s-2\int_0^t b^* \tilde{a}^{-1} b(X_s)ds\right)\right] \nonumber\\
 & &\ \ \ \  \cdot E_x\left[\exp\left(\int_0^t b^* \tilde{a}^{-1} b(X_s)ds + 2\int_0^t c(X_s)ds \right)g^2(X_t);
t<\tau_D\right]dx \nonumber\\
 &=&\int_{D}g^2(x)\hat E_x\left[\exp\left(\int_0^t b^* \tilde{a}^{-1} b(X_s)ds + 2\int_0^t c(X_s)ds \right);
t<\tau_D\right]dx  \nonumber\\
 &\leq&\varrho_2e^{\varrho_2t}\int_{D}g^2(x)dx
 \end{eqnarray}
where $\varrho_2>0$ is a constant independent of $g$. This gives
the existence of the extension of $P_t$ to $L^2(D;dx)$. Since
$C_b(D)$ is dense in $L^2(D;dx)$ and for $g\in C_b(D),
P_tg(x)\rightarrow g(x)$ as $t\rightarrow0$, the continuity
property of $P_t$ follows from (\ref{223}).

Define $$S_t=\exp\left(\int_0^t(\tilde
a^{-1}b)^*(X_s)dM_s-\frac{1}{2}\int_0^t b^* \tilde{a}^{-1}
b(X_s)ds + \int_0^t c(X_s)ds\right)$$ and $$\bar
M_t=\int_0^t(\tilde a^{-1}b)^*(X_s)dM_s.$$ Then
$S_t=1+\int_0^tS_sd\bar M_s+\int_0^tS_sc(X_s)ds$. By Ito's
formula, we obtain that for $u\in D({\cal L})$ and $t<\tau_D$,
$$u(X_t)S_t=u(X_0)+\int_0^tS_sdM_s^u+\int_0^tu(X_s)S_sd\bar M_s+\int_0^tS_s{\cal L}u(X_s)ds.$$
Following the argument of the proof of \cite[Theorem 3.2]{Zhang1},
we can show that $\{P_t\}_{t\geq0}$ coincides with
$\{T_t\}_{t\geq0}$ for this case. \vskip 0.2cm
 \noindent {\bf Case
2: $\hat b\in C_0^\infty(D)$}.

Similar to the proof of \cite[Theorem 3.3]{Zhang1}, we can show
that for $g\in L^2(D;dx)$,
\begin{eqnarray*}
& & T_tg(x)\nonumber\\
&=&E_x\left[\exp\left(\int_0^t(\tilde a^{-1}b)^*(X_s)dM_s-\frac{1}{2}\int_0^tb^*\tilde a^{-1}b(X_s)ds\right.\right.\nonumber\\
& &\ \ \ \ \left.\left. +\int_0^t c(X_s)ds-\int_0^t {\rm
div}\,\hat b(X_s)ds\right)g(X_t); t<\tau_D\right].
\end{eqnarray*}
The proof of this case is complete by (\ref{2004}).

 \vskip 0.2cm
 \noindent {\bf Case
3: $|\hat b|^2\in L^{p\vee 1}(D;dx)$}.

By Lemma \ref{lem111}(ii), we may choose a sequence $\{\hat b_n\in
C_0^\infty({\mathbb{R}^d})\}$ such that $|\hat b_n-\hat
b|^2\rightarrow 0 $ in $L^{p\vee 1}({\mathbb{R}^d};dx)$ and $\hat
b_n^H\rightarrow \hat b^H$ in $H^{1,2}({\mathbb{R}^d})$  as
$n\rightarrow\infty$.

Let $\{T_t^n\}_{t\ge 0}$ be the semigroup corresponding to the
quadratic form $\cal E$ with $\hat b_n$ in place of $\hat b$.
Then, for $f,g \in L^2(D;dx)$, we have
\begin{eqnarray} \label{29}
&&\int_{D} f(x)T_t^ng(x)dx\nonumber\\
&=&E_m\left[f(X_0)g(X_t)\exp\left(\int_0^t(\tilde a^{-1}b)^*(X_s)dM_s-\frac{1}{2}\int_0^tb^*\tilde{a}^{-1}b(X_s)ds\right.\right.\nonumber\\
& &\ \ \ \ \left.\left. +\int_0^t c(X_s)ds+N^{{\hat
b_n}^H}_t-\int_0^t {\hat b_n}^H(X_s)ds\right); t<\tau_D\right].
\end{eqnarray}
By \cite[Theorem 1.3]{Zhang2}, the left-hand side of (\ref{29})
converges to $\int_{D} f(x)T_tg(x)dx$ as $n\rightarrow\infty$.

We will prove below that the right-hand side  of (\ref{29})
converges to the right-hand side of (\ref{28}) as
$n\rightarrow\infty$. Define for $t\ge 0$,
\begin{eqnarray*}
Y^n_t&=&g(X_t)\exp\left(\int_0^t(\tilde a^{-1}b)^*(X_s)dM_s-\frac{1}{2}\int_0^tb^*\tilde{a}^{-1}b(X_s)ds\right.\nonumber\\
& &\ \ \ \ \left. +\int_0^t c(X_s)ds+N^{{\hat b_n}^H}_t-\int_0^t
{\hat b}_n^H(X_s)ds\right),\ \ n\in \mathbb{N},
\end{eqnarray*}
and
\begin{eqnarray*}
Y_t&=&g(X_t)\exp\left(\int_0^t(\tilde a^{-1}b)^*(X_s)dM_s-\frac{1}{2}\int_0^tb^*\tilde{a}^{-1}b(X_s)ds\right.\nonumber\\
& &\ \ \ \ \left. +\int_0^t c(X_s)ds+N^{{\hat b}^H}_t-\int_0^t
{\hat b}^H(X_s)ds\right).
\end{eqnarray*}

Then, the right-hand sides of (\ref{29}) and (\ref{28}) equal
$E_{f\cdot m}[Y^n_t; t<\tau_D]$ and $E_{f\cdot m}[Y_t; t<\tau_D]$,
respectively. To complete the proof, we need only show that
$\{Y^n_t1_{ t<\tau_D}\}$ is $P_{f\cdot m}$-uniformly integrable.
We will establish this below by proving that $\sup_{n\in
\mathbb{N}}E_{f\cdot m}[(Y^n_t)^2; t<\tau_D]<\infty$.

In fact, we obtain by Cauchy-Schwarz inequality that
\begin{eqnarray*}
& &E_{f\cdot m}[(Y^n_t)^2; t<\tau_D]\nonumber\\&=&E_{f\cdot
m}\left[g^2(X_t)\exp\left(2\int_0^t(\tilde
a^{-1}b)^*(X_s)dM_s-\int_0^tb^*\tilde{a}^{-1}b(X_s)ds\right.\right.\nonumber\\
& &\ \ \ \ \left.\left. +2\int_0^t c(X_s)ds+2N^{{\hat
b_n}^H}_t-2\int_0^t {\hat b}_n^H(X_s)ds\right);
t<\tau_D\right]\\
&=&E_{f\cdot m}\left[g^2(X_t)\exp\left(\frac{1}{2}\int_0^t(\tilde
a^{-1}b)^*(X_s)dM_s
-\frac{1}{4}\int_0^tb^*\tilde{a}^{-1}b(X_s)ds\right.\right. \nonumber\\
& &\ \ \ \ \left.\left. +\frac{1}{2}\int_0^t c(X_s)ds+2N^{{\hat b_n}^H}_t-2\int_0^t {\hat b}_n^H(X_s)ds\right) \right.\\
& &\left.\ \ \ \ \cdot \exp\left(\frac{3}{2}\int_0^t(\tilde
a^{-1}b)^*(X_s)dM_s-\frac{3}{4}\int_0^tb^*\tilde{a}^{-1}b(X_s)ds\right.\right.\\
& &\ \ \ \ \ \ \ \left.\left.+\frac{3}{2}\int_0^t c(X_s)ds\right);
t<\tau_D\right]\\
&\leq&E_{f\cdot m}\left[g^4(X_t)\exp\left(\int_0^t(\tilde
a^{-1}b)^*(X_s)dM_s
-\frac{1}{2}\int_0^tb^*\tilde{a}^{-1}b(X_s)ds\right.\right. \nonumber\\
& &\ \ \ \ \left.\left. +\int_0^t c(X_s)ds+N^{{4\hat
b_n}^H}_t-\int_0^t 4{\hat b}_n^H(X_s)ds\right);
t<\tau_D\right]^{1/2} \\
& &\cdot E_{f\cdot m}\left[\exp\left(3\int_0^t(\tilde
a^{-1}b)^*(X_s)dM_s\right.\right.\nonumber\\
& &\ \ \ \ \left.\left.
-\frac{3}{2}\int_0^tb^*\tilde{a}^{-1}b(X_s)ds +3\int_0^t
c(X_s)ds\right);
t<\tau_D\right]^{1/2}\\
&=&\left(\int_{D}f(x)T_t^{n'}g^4(x)dx\right)^{1/2} \\
& &\cdot E_{f\cdot m}\left[\exp\left(3\int_0^t(\tilde
a^{-1}b)^*(X_s)dM_s\right.\right.\nonumber\\
& &\ \ \ \ \left.\left.
-\frac{3}{2}\int_0^tb^*\tilde{a}^{-1}b(X_s)ds +3\int_0^t
c(X_s)ds\right); t<\tau_D\right]^{1/2},
\end{eqnarray*}
where $\{T_t^{n'}\}_{t\ge 0}$ is the semigroup corresponding to
the quadratic form ${\cal E}$ with $4\hat b_n$ in place of $\hat
b$. Thus, we obtain by \cite[Theorem 1.3]{Zhang2} and Lemma
\ref{432} that
\begin{eqnarray*}
& &\sup_{n\in \mathbb{N}} E_{f\cdot m}[(Y^n_t)^2;
t<\tau_D]\\
&\leq& \sup_{n\in \mathbb{N}}\left(\int_{D}f(x)T_t^{n'}g^4(x)dx\right)^{1/2}\cdot E_{f\cdot m}\left[\exp\left(3\int_0^t(\tilde a^{-1}b)^*(X_s)dM_s \right.\right. \\
& &\left.\left. \ \ \ \ \ \ \ \
-\frac{3}{2}\int_0^tb^*\tilde{a}^{-1}b(X_s)ds+3\int_0^t
c(X_s)ds\right); t<\tau_D\right]^{1/2}\\
&<&\infty.
\end{eqnarray*}
\hfill\fbox

\bigskip

{ \noindent {\bf\large Acknowledgments} \vskip 0.1cm  \noindent We
acknowledge the support of NSFC (Grant No. 11361021) and NSERC
(Grant No. 311945-2013). We thank the anonymous referee and Professors Z.Q. Chen, K. Kuwae, T.S. Zhang for
very helpful comments, which improved the presentation of the paper.}


\begin{thebibliography}{1234}

\bibitem{Aron1} D.G. Aronson,  Bounds on the fundamental solution of a parabolic equation, Bull. Amer. Math. Soc. 73 (1967) 890-896.

    \bibitem{Aron} D.G. Aronson,  Non-negative solutions of linear parabolic equations, Ann. Scuola Norm. Sup. Pisa 22 (1968) 607-694.


\bibitem{CMS} C.Z. Chen, L. Ma, W. Sun,  Stochastic calculus for Markov processes associated
with non-symmetric Dirichlet forms, Sci. China Math. { 55} (2012)
2195-2203.

\bibitem{Zhang3} Z.Q. Chen, T.S. Zhang,  Time-reversal and elliptic boundary value problems, Ann. Prob. { 37} (2009)
1008-1043.

\bibitem{Zhang5} Z.Q. Chen, T.S. Zhang,  A probabilistic approach to mixed boundary value problems for elliptic operators with
singular coefficients, Proc. Amer. Math. Soc. {142} (2014)
2135-2149.

\bibitem{Chen} Z.Q.
Chen, Z. Zhao, Diffusion processes and second order elliptic
operators with singular coefficients for lower order terms, Math.
Ann.  { 302} (1995) 323-357.

\bibitem{Cho} S. Cho, P. Kim, H. Park, Two-sided estimates on Dirichlet heat
kernels for time-dependent parabolic operators with singular
drifts in $C^{1,\alpha}$-domains, J. Diff. Equat. { 252} (2012)
1101-1145.

\bibitem{Faz}G. Di Fazio, $L^p$ estimates for divergence form elliptic equations with discontinuous
coeffcients, Boll. Un. Mat. Ital. A(7) 10 (1996) 409-420.

\bibitem{Fu94} {M. Fukushima, Y. Oshima, M. Takeda}, {{Dirichlet Forms
              and Symmetric Markov Processes}}, Second revised and extended edition, Walter de
              Gruyter,
            2011.

\bibitem{Kakutani} S. Kakutani, Two-dimensional Brownian motion and harmonic functions,
 Proc. Imp. Acad. Tokyo {20} (1944) 706-714.

\bibitem{Kato} T. Kato, Perturbation Theory for Linear Operators, Springer-Verlag,
1980.

\bibitem{Kenig} C. Kenig, H. Koch, J. Pipher, T. Toro, A new approach to absolute
continuity of elliptic measure, with applications to non-symmetric
equations, Adv. Math. { 153} (2000) 231-298.

\bibitem{KS} P. Kim, R.M. Song, Two-sided estimates on the density of Brownian motion with singular drift, Ill. J. Math. {
50} (2006) 635-688.

\bibitem{Lierl} J. Lierl, L. Saloff-Coste, Parabolic Harnack inequality for time-dependent non-symmetric Dirichlet
forms, arXiv:1205.6493v4.

\bibitem{Zhang1} J. Lunt, T.J. Lyons, T.S. Zhang,  Integrability of functonals of Dirichlet processes,
probabilistic representations of semigroups, and estimates of heat
kernels, J. Funct. Anal. { 153}  (1998) 320-342.

\bibitem{MR92} Z.M. Ma, M. R\"ockner, {Introduction to the Theory
              of (Non-symmetric) Dirichlet Forms},
              Springer-Verlag,       1992.

\bibitem{Meyers} N.G. Meyers, An $L^p$-estimate for the gradient of solutions of second order elliptic divergence
equations, Ann. Scuola Norm. Sup. Pisa {17} (1963) 189-206.

\bibitem{Morrey} C.B.Jr. Morrey, Second order elliptic equations in several variables and H\"older continuity, Math Z.
{72} (1959) 146-164.

\bibitem{M2} C.B.Jr. Morrey, Multiple
Integrals in the Calculus of Variations, Springer-Verlag, 1966.

\bibitem{Moser} J. Moser, A Harnack inequality for parabolic differential equations,
 Comm. Pure Appl. Math. { 17} (1964) 101-134.

\bibitem{N} S. Nakao, Stochastic calculus for continuous additive functionals of zero
energy,
 { Z. Wahrsch. verw. Gebiete} { 68} (1985) 557-578.

\bibitem{oshima} {Y. Oshima}, {Lecture on Dirichlet Spaces}, Univ. Erlangen-N\"urnberg, 1988.

\bibitem{Zhang2} M. R\"ockner, T.S. Zhang, Convergence of operator semigroups generated by elliptic
operators, Osaka J. Math. { 34} (1997) 923-932.

\bibitem{Trudinger} N. Trudinger, Linear elliptic operators
with measurable coefficients, Ann. Scuola Norm. Sup. Pisa { 27}
(1973) 255-308.

\bibitem{Alexander 2013} A. Walsh,
Stochastic integration with respect to additive functionals of
zero quadratic variation, { Bernoulli} { 19} (2013) 2414-2436.

\bibitem{Zhang4} T.S. Zhang,  A probabilistic approach to
Dirichlet problems of semilinear elliptic PDEs with singular
coefficients,  Ann. Probab. { 39} (2011) 1502-1527.

\end{thebibliography}
\end{document}